\documentclass[reqno,a4paper]{amsart}
\usepackage{amssymb}
\usepackage[colorlinks=true,hypertex]{hyperref}%,pagebackref
\allowdisplaybreaks[4]
\theoremstyle{plain}
\newtheorem{thm}{Theorem}[section]

\newtheorem{prop}{Proposition}[section]
\theoremstyle{remark}
\newtheorem{rem}{Remark}[subsection]
\DeclareMathOperator{\td}{d\mspace{-2mu}}
\DeclareMathOperator{\re}{Re}
\date{Drafted on 7 April 2008 in Melbourne}
\numberwithin{equation}{section}
\date{Completed on Sunday 22 September 2008 in Carlisle B, VU Student Village, Australia}
\date{}

\begin{document}

\title[Bounds for the ratio of two gamma functions]
{Bounds for the ratio of two gamma functions---From Wendel's and related inequalities to logarithmically completely monotonic functions}

\author[F. Qi]{Feng Qi}
\address[F. Qi]{Research Institute of Mathematical Inequality Theory, Henan Polytechnic University, Jiaozuo City, Henan Province, 454010, China}
\email{\href{mailto: F. Qi <qifeng618@gmail.com>}{qifeng618@gmail.com}, \href{mailto: F. Qi <qifeng618@hotmail.com>}{qifeng618@hotmail.com}, \href{mailto: F. Qi <qifeng618@qq.com>}{qifeng618@qq.com}}
\urladdr{\url{http://qifeng618.spaces.live.com}}

\begin{abstract}
In this expository and survey paper, along one of main lines of bounding the ratio of two gamma functions, we look back and analyse some inequalities, several complete monotonicity of functions involving ratios of two gamma or $q$-gamma functions, and necessary and sufficient conditions for functions involving ratios of two gamma or $q$-gamma functions to be logarithmically completely monotonic.
\end{abstract}



\thanks{The author was partially supported by the China Scholarship Council}

\thanks{This paper was typeset using \AmS-\LaTeX}

\maketitle

\tableofcontents

\section{Introduction}

\subsection{The gamma and $q$-gamma functions}
It is well-known that the classical Euler gamma function may be defined by
\begin{equation}\label{egamma}
\Gamma(x)=\int^\infty_0t^{x-1} e^{-t}\td t,\quad x>0.
\end{equation}
The logarithmic derivative of $\Gamma(x)$, denoted by $\psi(x)=\frac{\Gamma'(x)}{\Gamma(x)}$, is called the psi or digamma function, and $\psi^{(k)}(x)$ for $k\in \mathbb{N}$ are called the polygamma functions. It is common knowledge that special functions $\Gamma(x)$, $\psi(x)$ and $\psi^{(k)}(x)$ for $k\in\mathbb{N}$ are fundamental and important and have much extensive applications in mathematical sciences.
\par
The $q$-analogues of $\Gamma$ and $\psi$ are defined \cite[pp.~493\nobreakdash--496]{andrews} for $x>0$ by
\begin{gather}\label{q-gamma-dfn}
\Gamma_q(x)=(1-q)^{1-x}\prod_{i=0}^\infty\frac{1-q^{i+1}}{1-q^{i+x}},\quad 0<q<1,\\
\label{q-gamma-dfn-q>1}
\Gamma_q(x)=(q-1)^{1-x}q^{\binom{x}2}\prod_{i=0}^\infty\frac{1-q^{-(i+1)}}{1-q^{-(i+x)}}, \quad q>1,
\end{gather}
and
\begin{align}\label{q-gamma-1.4}
\psi_q(x)=\frac{\Gamma_q'(x)}{\Gamma_q(x)}&=-\ln(1-q)+\ln q \sum_{k=0}^\infty\frac{q^{k+x}}{1-q^{k+x}}\\
&=-\ln(1-q)-\int_0^\infty\frac{e^{-xt}}{1-e^{-t}}\td\gamma_q(t) \label{q-gamma-1.5}
\end{align}
for $0<q<1$, where $\td\gamma_q(t)$ is a discrete measure with positive masses $-\ln q$ at the positive points $-k\ln q$ for $k\in\mathbb{N}$, more accurately,
\begin{equation}
\gamma_q(t)=
\begin{cases}
-\ln q\sum\limits_{k=1}^\infty\delta(t+k\ln q),&0<q<1,\\ t,&q=1.
\end{cases}
\end{equation}
See~\cite[p.~311]{Ismail-Muldoon-119}.
\par
The $q$-gamma function $\Gamma_q(z)$ has the following basic properties:
\begin{equation}
\lim_{q\to1^+}\Gamma_q(z)=\lim_{q\to1^-}\Gamma_q(z)=\Gamma(z)\quad \text{and}\quad \Gamma_q(x)=q^{\binom{x-1}2}\Gamma_{1/q}(x).
\end{equation}

\subsection{The definition and properties of completely monotonic functions}
A function $f$ is said to be completely monotonic on an interval $I$ if $f$ has derivatives of all orders on $I$ and $(-1)^{n}f^{(n)}(x)\ge0$ for $x \in I$ and $n \ge0$.
\par
The class of completely monotonic functions has the following basic properties.

\begin{thm}[{\cite[p.~161]{widder}}]\label{p.161-widder}
A necessary and sufficient condition that $f(x)$ should be completely monotonic for $0<x<\infty$ is that
\begin{equation}
f(x)=\int_0^\infty e^{-xt}\td\alpha(t),
\end{equation}
where $\alpha(t)$ is nondecreasing and the integral converges for $0<x<\infty$.
\end{thm}

\begin{thm}[{\cite[p.~83]{bochner}}]\label{p.83-bochner}
If $f(x)$ is completely monotonic on $I$, $g(x)\in I$, and $g'(x)$ is completely monotonic on $(0,\infty)$, then $f(g(x))$ is completely monotonic on $(0,\infty)$.
\end{thm}

\subsection{The logarithmically completely monotonic functions}
A positive and $k$-times differentiable function $f(x)$ is said to be $k$-log-convex (or $k$-log-concave, respectively) with $k\ge2$ on an interval $I$ if and only if $[\ln f(x)]^{(k)}$ exists and $[\ln f(x)]^{(k)}\ge0$ (or $[\ln f(x)]^{(k)}\le0$, respectively) on $I$.
\par
A positive function $f(x)$ is said to be logarithmically completely monotonic on an interval $I\subseteq\mathbb{R}$ if it has derivatives of all orders on $I$ and its logarithm $\ln f(x)$ satisfies $(-1)^k[\ln f(x)]^{(k)}\ge0$ for $k\in\mathbb{N}$ on $I$.
\par
The notion ``logarithmically completely monotonic function'' was first put forward in~\cite{Atanassov} without an explicit definition. This terminology was explicitly recovered in~\cite{minus-one} whose revised and expanded version was formally published as~\cite{minus-one.tex-rev}.
\par
It has been proved once and again in \cite{CBerg, clark-ismail-NFAA.tex, clark-ismail.tex, compmon2, absolute-mon.tex, minus-one, minus-one.tex-rev, schur-complete} that a logarithmically completely monotonic function on an interval $I$ must also be completely monotonic on $I$. C.~Berg points out in~\cite{CBerg} that these functions are the same as those studied by Horn \cite{horn} under the name infinitely divisible completely monotonic functions. For more information, please refer to \cite{CBerg, auscmrgmia} and related references therein.

\subsection{Outline of this paper}
The history of bounding the ratio of two gamma functions has been longer than sixty years since the paper~\cite{wendel} by J. G. Wendel was published in 1948.
\par
The motivations of bounding the ratio of two gamma functions are various, including establishment of asymptotic relation, refinements of Wallis' formula, approximation to $\pi$, needs in statistics and other mathematical sciences.
\par
In this expository and survey paper, along one of main lines of bounding the ratio of two gamma functions, we look back and analyse some inequalities such as Wendel's double inequality, Kazarinoff's refinement of Wallis' formula, Watson's monotonicity, Gautschi's double inequality, and Kershaw's first double inequality, the complete monotonicity of several functions involving ratios of two gamma or $q$-gamma functions by Bustoz, Ismail, Lorch and Muldoon, and necessary and sufficient conditions for functions involving ratios of two gamma or $q$-gamma functions to be logarithmically completely monotonic.

\section{Some inequalities for bounding the ratio of two gamma functions}

In this section, we look back and analyse some related inequalities for bounding the ratio of two gamma functions.

\subsection{Wendel's double inequality}
Our starting point is a paper published in 1948 by J. G. Wendel, which is the earliest one we can search out to the best of our ability.
\par
In order to establish the classical asymptotic relation
\begin{equation}\label{wendel-approx}
\lim_{x\to\infty}\frac{\Gamma(x+s)}{x^s\Gamma(x)}=1
\end{equation}
for real $s$ and $x$, by using H\"older's inequality for integrals, J. G. Wendel~\cite{wendel} proved elegantly the double inequality
\begin{equation}\label{wendel-inequal}
\biggl(\frac{x}{x+s}\biggr)^{1-s}\le\frac{\Gamma(x+s)}{x^s\Gamma(x)}\le1
\end{equation}
for $0<s<1$ and $x>0$.

\begin{rem}\label{rem-2.1.1}
The inequality~\eqref{wendel-inequal} can be rewritten for $0<s<1$ and $x>0$ as
\begin{equation}\label{wendel-inequal-rew}
(x+s)^{s-1}\frac{\Gamma(x+1)}{\Gamma(x+s)}\le1\le x^{s-1}\frac{\Gamma(x+1)}{\Gamma(x+s)}.
\end{equation}
The relation~\eqref{wendel-approx} results in
\begin{equation}
\lim_{x\to\infty}(x+s)^{s-1}\frac{\Gamma(x+1)}{\Gamma(x+s)} =\lim_{x\to\infty}x^{s-1}\frac{\Gamma(x+1)}{\Gamma(x+s)}=1
\end{equation}
which hints us that the functions
\begin{equation}\label{2.5-function}
(x+s)^{s-1}\frac{\Gamma(x+1)}{\Gamma(x+s)}\quad \text{and}\quad x^{s-1}\frac{\Gamma(x+1)}{\Gamma(x+s)}
\end{equation}
or
\begin{equation}\label{2.5-function-2}
(x+s)\biggl[\frac{\Gamma(x+1)}{\Gamma(x+s)}\biggr]^{1/(s-1)}\quad \text{and}\quad x\biggl[\frac{\Gamma(x+1)}{\Gamma(x+s)}\biggr]^{1/(s-1)}
\end{equation}
are possibly increasing and decreasing respectively.
\end{rem}

\begin{rem}
In~\cite[p.~257, 6.1.46]{abram}, the following limit was listed: For real numbers $a$ and $b$,
\begin{equation}\label{gamma-ratio-lim}
\lim_{x\to\infty}\biggl[x^{b-a}\frac{\Gamma(x+a)}{\Gamma(x+b)}\biggr]=1.
\end{equation}
The limits~\eqref{wendel-approx} and~\eqref{gamma-ratio-lim} are equivalent to each other since
\begin{equation*}
x^{t-s}\frac{\Gamma(x+s)}{\Gamma(x+t)}=\frac{\Gamma(x+s)}{x^s\Gamma(x)} \cdot\frac{x^t\Gamma(x)}{\Gamma(x+t)}.
\end{equation*}
Hence, the limit~\eqref{gamma-ratio-lim} is presumedly called as Wendel's limit.
\end{rem}

\begin{rem}
Due to unknown reasons, Wendel's paper~\cite{wendel} was seemingly neglected by nearly all mathematicians for more than fifty years until it was mentioned in~\cite{Merkle-JMAA-99}, to the best of my knowledge.
\end{rem}

\subsection{Kazarinoff's double inequality}\label{Wallis-section}

Starting from
\begin{equation}\label{John-Wallis-ineq}
\frac1{\sqrt{\pi(n+1/2)}}<\frac{(2n-1)!!}{(2n)!!}<\frac1{\sqrt{\pi n}},\quad n\in\mathbb{N},
\end{equation}
one form of the celebrated formula of John Wallis, which had been quoted for more than a century before 1950s by writers of textbooks, D. K. Kazarinoff proved in~\cite{Kazarinoff-56} that the sequence $\theta(n)$ defined by
\begin{equation}\label{theta-dfn-kazar}
\frac{(2n-1)!!}{(2n)!!} =\frac1{\sqrt{\pi[n+\theta(n)]}}
\end{equation}
satisfies $\frac14<\theta(n)<\frac12$ for $n\in\mathbb{N}$, that is,
\begin{equation}\label{Wallis'inequality}
\frac1{\sqrt{\pi(n+1/2)}}<\frac{(2n-1)!!}{(2n)!!}
<\frac1{\sqrt{\pi(n+1/4)}},\quad n\in\mathbb{N}.
\end{equation}

\begin{rem}
It was said in~\cite{Kazarinoff-56} that it is unquestionable that inequalities similar to~\eqref{Wallis'inequality} can be improved indefinitely but at a sacrifice of simplicity, which is why the inequality~\eqref{John-Wallis-ineq} had survived so long.
\end{rem}

\begin{rem}
Kazarinoff's proof of~\eqref{Wallis'inequality} is based upon the property
\begin{equation}\label{Phi-ineq}
[\ln\phi(t)]''-\{[\ln\phi(t)]'\}^2>0
\end{equation}
of the function
\begin{equation}
\phi(t)=\int_0^{\pi/2}\sin^tx\td x=\frac{\sqrt\pi\,}2\cdot\frac{\Gamma((t+1)/2)}{\Gamma((t+2)/2)}
\end{equation}
for $-1<t<\infty$. The inequality~\eqref{Phi-ineq} was proved by making use of the well-known Legendre's formula
\begin{equation}\label{Legendre's-formula}
\psi(x)=-\gamma+\int_0^1\frac{t^{x-1}-1}{t-1}\td t
\end{equation}
for $x>0$ and estimating the integrals
\begin{equation}
\int_0^1\frac{x^t}{1+x}\td x\quad\text{and}\quad \int_0^1\frac{x^t\ln x}{1+x}\td x.
\end{equation}
Since~\eqref{Phi-ineq} is equivalent to the statement that the reciprocal of $\phi(t)$ has an everywhere negative second derivative, therefore, for any positive $t$, $\phi(t)$ is less than the harmonic mean of $\phi(t-1)$ and $\phi(t+1)$, which implies
\begin{equation}\label{karz-2.17-ineq}
\frac{\Gamma((t+1)/2)}{\Gamma((t+2)/2)}<\frac2{\sqrt{2t+1}},\quad t>-\frac12.
\end{equation}
As a subcase of this result, the right-hand side inequality in~\eqref{Wallis'inequality} is established.
\end{rem}

\begin{rem}
Replacing $t$ by $2t$ in~\eqref{karz-2.17-ineq} and rearranging yield
\begin{equation}\label{karz-2.17-ineq-rew}
\frac{\Gamma(t+1)}{\Gamma(t+1/2)}>\sqrt{t+\frac14}\quad\Longleftrightarrow\quad \biggl(t+\frac14\biggr)^{1/2-1}\frac{\Gamma(t+1)}{\Gamma(t+1/2)}>1
\end{equation}
for $t>-\frac14$. From~\eqref{wendel-approx}, it follows that
\begin{equation}
\lim_{x\to\infty}\biggl(t+\frac14\biggr)^{1/2-1}\frac{\Gamma(t+1)}{\Gamma(t+1/2)}=1.
\end{equation}
This suggests that the function
\begin{equation}
\biggl(t+\frac14\biggr)^{1/2-1}\frac{\Gamma(t+1)}{\Gamma(t+1/2)}\quad \text{or}\quad \biggl(t+\frac14\biggr)\biggl[\frac{\Gamma(t+1)}{\Gamma(t+1/2)}\biggr]^{1/(1/2-1)}
\end{equation}
is perhaps decreasing, more strongly, logarithmically completely monotonic.
\end{rem}

\begin{rem}\label{kazarinoff-rem-5}
The inequality~\eqref{Phi-ineq} may be rewritten as
\begin{equation}
\psi'\biggl(\frac{t+1}2\biggr)-\psi'\biggl(\frac{t+2}2\biggr) >\biggl[\psi\biggl(\frac{t+1}2\biggr)-\psi\biggl(\frac{t+2}2\biggr)\biggr]^2
\end{equation}
for $t>-1$. Letting $u=\frac{t+1}2$ in the above inequality yields
\begin{equation}
\psi'(u)-\psi'\biggl(u+\frac12\biggr) >\biggl[\psi(u)-\psi\biggl(u+\frac12\biggr)\biggr]^2
\end{equation}
for $u>0$. This inequality has been generalized in~\cite{Comp-Mon-Digamma-Trigamma-Divided.tex} to the complete monotonicity of a function involving divided differences of the digamma and trigamma functions as follows.

\begin{thm}\label{CMDT-divided-thm}
For real numbers $s$, $t$, $\alpha=\min\{s,t\}$ and $\lambda$, let
\begin{equation}\label{Delta-lambda-dfn}
\Delta_{s,t;\lambda}(x)=\begin{cases}\bigg[\dfrac{\psi(x+t) -\psi(x+s)}{t-s}\bigg]^2
+\lambda\dfrac{\psi'(x+t)-\psi'(x+s)}{t-s},&s\ne t\\
[\psi'(x+s)]^2+\lambda\psi''(x+s),&s=t
\end{cases}
\end{equation}
on $(-\alpha,\infty)$. Then the function $\Delta_{s,t;\lambda}(x)$ has the following complete monotonicity:
\begin{enumerate}
\item
For $0<|t-s|<1$,
\begin{enumerate}
\item
the function $\Delta_{s,t;\lambda}(x)$ is completely monotonic on $(-\alpha,\infty)$ if and only if $\lambda\le1$,
\item
so is the function $-\Delta_{s,t;\lambda}(x)$ if and only if $\lambda\ge\frac1{|t-s|}$;
\end{enumerate}
\item
For $|t-s|>1$,
\begin{enumerate}
\item
the function $\Delta_{s,t;\lambda}(x)$ is completely monotonic on $(-\alpha,\infty)$ if and only if $\lambda\le\frac1{|t-s|}$,
\item
so is the function $-\Delta_{s,t;\lambda}(x)$ if and only if $\lambda\ge1$;
\end{enumerate}
\item
For $s=t$, the function $\Delta_{s,s;\lambda}(x)$ is completely monotonic on $(-s,\infty)$ if and only if $\lambda\le1$;
\item
For $|t-s|=1$,
\begin{enumerate}
\item
the function $\Delta_{s,t;\lambda}(x)$ is completely monotonic if and only if $\lambda<1$,
\item
so is the function $-\Delta_{s,t;\lambda}(x)$ if and only if $\lambda>1$,
\item
and $\Delta_{s,t;1}(x)\equiv0$.
\end{enumerate}
\end{enumerate}
\end{thm}
Taking $\lambda=s-t>0$ in Theorem~\ref{CMDT-divided-thm} produces that the function $\frac{\Gamma(x+s)}{\Gamma(x+t)}$ on $(-t,\infty)$ is increasingly convex if $s-t>1$ and increasingly concave if $0<s-t<1$.
\end{rem}

\subsection{Watson's monotonicity}\label{Watson-sec}

In 1959, motivated by the result in~\cite{Kazarinoff-56} mentioned in Section~\ref{Wallis-section}, G.~N.~Watson~\cite{waston} observed that
\begin{multline}\label{watson-formula}
\frac1x\cdot\frac{[\Gamma(x+1)]^2}{[\Gamma(x+1/2)]^2} ={}_2F_1\biggl(-\frac12,-\frac12;x;1\biggr)\\*
=1+\frac1{4x}+\frac1{32x(x+1)} +\sum_{r=3}^\infty\frac{[(-1/2)\cdot(1/2)\cdot(3/2)\cdot(r-3/2)]^2} {r!x(x+1)\dotsm(x+r-1)}
\end{multline}
for $x>-\frac12$, which implies the much general function
\begin{equation}\label{theta-dfn}
\theta(x)=\biggl[\frac{\Gamma(x+1)}{\Gamma(x+1/2)}\biggr]^2-x
\end{equation}
for $x>-\frac12$, whose special case is the sequence $\theta(n)$ for $n\in\mathbb{N}$ defined in~\eqref{theta-dfn-kazar}, is decreasing and
\begin{equation}
\lim_{x\to\infty}\theta(x)=\frac14\quad \text{and}\quad \lim_{x\to(-1/2)^+}\theta(x)=\frac12.
\end{equation}
This implies apparently the sharp inequalities
\begin{equation}\label{theta-l-u-b}
\frac14<\theta(x)<\frac12
\end{equation}
for $x>-\frac12$,
\begin{equation}\label{watson-special-ineq}
\sqrt{x+\frac14}\,< \frac{\Gamma(x+1)}{\Gamma(x+1/2)}\le \sqrt{x+\frac14+\biggl[\frac{\Gamma(3/4)}{\Gamma(1/4)}\biggr]^2}\, =\sqrt{x+0.36423\dotsm}
\end{equation}
for $x\ge-\frac14$, and, by Wallis cosine formula~\cite{WallisFormula.html},
\begin{equation}\label{best-bounds-Wallis}
\frac{1}{\sqrt{\pi(n+{4}/{\pi}-1)}}\le\frac{(2n-1)!!}{(2n)!!}
<\frac{1}{\sqrt{\pi(n+1/4)}},\quad n\in\mathbb{N}.
\end{equation}

\begin{rem}
In~\cite{waston}, an alternative proof of the double inequality~\eqref{theta-l-u-b} was also provided.
\end{rem}

\begin{rem}
It is easy to see that the inequality~\eqref{watson-special-ineq} extends and improves~\eqref{wendel-inequal} when $s=\frac12$.
\end{rem}

\begin{rem}
The left-hand side inequality in~\eqref{best-bounds-Wallis} is better than the corresponding one in~\eqref{Wallis'inequality}.
\end{rem}

\begin{rem}
The formula~\eqref{watson-formula} implies the complete monotonicity of the function $\theta(x)$ on $\bigl(-\frac12,\infty\bigr)$ defined by~\eqref{theta-dfn}.
\end{rem}

\subsection{Gautschi's double inequalities}

The first result of the paper~\cite{gaut} was the double inequality
\begin{equation}\label{gaut-3-ineq}
\frac{(x^p+2)^{1/p}-x}2<e^{x^p}\int_x^\infty e^{-t^p}\td t\le c_p\biggl[\biggl(x^p+\frac1{c_p}\biggr)^{1/p}-x\biggr]
\end{equation}
for $x\ge0$ and $p>1$, where
\begin{equation}
c_p=\biggl[\Gamma\biggl(1+\frac1p\biggr)\biggr]^{p/(p-1)}
\end{equation}
or $c_p=1$. By an easy transformation, the inequality~\eqref{gaut-3-ineq} was written in terms of the complementary gamma function
\begin{equation}
\Gamma(a,x)=\int_x^\infty e^{-t}t^{a-1}\td t
\end{equation}
as
\begin{equation}\label{gaut-4-ineq}
\frac{p[(x+2)^{1/p}-x^{1/p}]}2<e^x\Gamma\biggl(\frac1p,x\biggr)\le pc_p\biggl[\biggl(x+\frac1{c_p}\biggr)^{1/p}-x^{1/p}\biggr]
\end{equation}
for $x\ge0$ and $p>1$. In particular, if letting $p\to\infty$, the double inequality
\begin{equation}
\frac12\ln\biggl(1+\frac2x\biggr)\le e^xE_1(x)\le\ln\biggl(1+\frac1x\biggr)
\end{equation}
for the exponential integral $E_1(x)=\Gamma(0,x)$ for $x>0$ was derived from~\eqref{gaut-4-ineq}, in which the bounds exhibit the logarithmic singularity of $E_1(x)$ at $x=0$. As a direct consequence of the inequality~\eqref{gaut-4-ineq} for $p=\frac1s$, $x=0$ and $c_p=1$, the following simple inequality for the gamma function was deduced:
\begin{equation}\label{gaut-none-ineq}
2^{s-1}\le\Gamma(1+s)\le1,\quad 0\le s\le 1.
\end{equation}
\par
The second result of the paper~\cite{gaut} was a sharper and more general inequality
\begin{equation}\label{gaut-6-ineq}
e^{(s-1)\psi(n+1)}\le\frac{\Gamma(n+s)}{\Gamma(n+1)}\le n^{s-1}
\end{equation}
for $0\le s\le1$ and $n\in\mathbb{N}$ than~\eqref{gaut-none-ineq}. It was obtained by proving that the function
\begin{equation}
f(s)=\frac1{1-s}\ln\frac{\Gamma(n+s)}{\Gamma(n+1)}
\end{equation}
is monotonically decreasing for $0\le s<1$. Since $\psi(n)<\ln n$, it was derived from the inequality~\eqref{gaut-6-ineq} that
\begin{equation}\label{gaut-6-ineq-simp}
\biggl(\frac1{n+1}\biggr)^{1-s}\le\frac{\Gamma(n+s)}{\Gamma(n+1)}\le\biggl(\frac1n\biggr)^{1-s}, \quad 0\le s\le1,
\end{equation}
which was also rewritten as
\begin{equation}\label{euler-gaut}
\frac{n!(n+1)^{s-1}}{(s+1)(s+2)\dotsm(s+n-1)}\le\Gamma(1+s) \le\frac{(n-1)!n^s}{(s+1)(s+2)\dotsm(s+n-1)},
\end{equation}
and so a simple proof of Euler's product formula in the segment $0\le s\le1$ was shown by letting $n\to\infty$ in~\eqref{euler-gaut}.

\begin{rem}
For more information on refining the inequality~\eqref{gaut-3-ineq}, please refer to~\cite{incom-gamma-L-N, qi-senlin-mia, Qi-Mei-99-gamma} and related references therein.
\end{rem}

\begin{rem}
The double inequalities~\eqref{gaut-6-ineq} and~\eqref{gaut-6-ineq-simp} can be rearranged as
\begin{equation}\label{gaut-ineq-1}
n^{1-s}\le\frac{\Gamma(n+1)}{\Gamma(n+s)}\le\exp((1-s)\psi(n+1))
\end{equation}
and
\begin{equation}\label{gaut-ineq-2}
n^{1-s}\le\frac{\Gamma(n+1)}{\Gamma(n+s)}\le (n+1)^{1-s}
\end{equation}
for $n\in\mathbb{N}$ and $0\le s\le 1$. Furthermore, the inequality~\eqref{gaut-ineq-2} can be rewritten as
\begin{equation}\label{gau-rew-1}
n^{1-s}\frac{\Gamma(n+s)}{\Gamma(n+1)}\le1\le (n+1)^{1-s}\frac{\Gamma(n+s)}{\Gamma(n+1)}
\end{equation}
or
\begin{equation}\label{gau-rew-2}
n\biggl[\frac{\Gamma(n+s)}{\Gamma(n+1)}\biggr]^{1/(1-s)}\le1\le (n+1)\biggl[\frac{\Gamma(n+s)}{\Gamma(n+1)}\biggr]^{1/(1-s)}.
\end{equation}
This supply us some possible clues to see that the sequences at the very ends of the inequalities~\eqref{gau-rew-1} and~\eqref{gau-rew-2} are monotonic.
\end{rem}

\begin{rem}
The left-hand side inequality in~\eqref{wendel-inequal} and the upper bound in~\eqref{gaut-ineq-1} have the following relationship
\begin{equation}\label{wendel-gautschi-comp}
(n+s)^{1-s}\le\exp((1-s)\psi(n+1))
\end{equation}
for $0\le s\le\frac12$ and $n\in\mathbb{N}$, and the inequality~\eqref{wendel-gautschi-comp} reverses for $s>e^{1-\gamma}-1=0.52620\dotsm$, since the function
\begin{equation}\label{Q(x)-dfn}
Q(x)=e^{\psi(x+1)}-x
\end{equation}
was proved in~\cite[Theorem~2]{Infinite-family-Digamma.tex} to be strictly decreasing on $(-1,\infty)$ and
\begin{equation}\label{Q-infty-lim}
\lim_{x\to\infty}Q(x)=\frac12.
\end{equation}
This means that Wendel's double inequality~\eqref{wendel-inequal} and Gautschi's first double inequality~\eqref{gaut-ineq-1} are not included each other, but they all contain Gautschi's second double inequality~\eqref{gaut-ineq-2}.
\end{rem}

\begin{rem}
In the reviews on the paper~\cite{gaut} by the Mathematical Reviews and the Zentralblatt MATH, there is no a word to comment on inequalities in~\eqref{gaut-ineq-1} and~\eqref{gaut-ineq-2}. However, these two double inequalities later became a major source of a large amount of study on bounding the ratio of two gamma functions.
\end{rem}

\subsection{Kershaw's first double inequality}\label{kershaw-sec}

Inspired by the inequality~\eqref{gaut-6-ineq-simp}, among other things, D. Kershaw presented in~\cite{kershaw} the following double inequality
\begin{equation}\label{gki1}
\biggl(x+\frac{s}2\biggr)^{1-s}<\frac{\Gamma(x+1)}{\Gamma(x+s)}
<\biggl[x-\frac12+\biggl(s+\frac14\biggr)^{1/2}\biggr]^{1-s}
\end{equation}
for $0<s<1$ and $x>0$. In the literature, it is called as Kershaw's first double inequality for the ratio of two gamma functions.

\begin{rem}
It is easy to see that the inequality~\eqref{gki1} refines and extends the inequality~\eqref{wendel-inequal} and~\eqref{gaut-ineq-2}.
\end{rem}

\begin{rem}
The inequality~\eqref{gki1} may be rearranged as
\begin{equation}\label{gki1-rew}
\biggl[x-\frac12+\biggl(s+\frac14\biggr)^{1/2}\biggr]^{s-1}\frac{\Gamma(x+1)}{\Gamma(x+s)} <1<\biggl(x+\frac{s}2\biggr)^{s-1}\frac{\Gamma(x+1)}{\Gamma(x+s)}
\end{equation}
for $x>0$ and $0<s<1$. By virtue of~\eqref{wendel-approx} or~\eqref{gamma-ratio-lim}, it is easy to see that
\begin{equation}
\lim_{x\to\infty}\biggl[x-\frac12+\biggl(s+\frac14\biggr)^{1/2}\biggr]^{s-1} \frac{\Gamma(x+1)}{\Gamma(x+s)} =\lim_{x\to\infty}\biggl(x+\frac{s}2\biggr)^{s-1}\frac{\Gamma(x+1)}{\Gamma(x+s)}=1.
\end{equation}
This insinuates the monotonicity, more strongly, the logarithmically complete monotonicity, of the functions
\begin{equation}
\biggl[x-\frac12+\biggl(s+\frac14\biggr)^{1/2}\biggr]^{s-1} \frac{\Gamma(x+1)}{\Gamma(x+s)} \quad\text{and}\quad \biggl(x+\frac{s}2\biggr)^{s-1}\frac{\Gamma(x+1)}{\Gamma(x+s)}
\end{equation}
or
\begin{equation}
\biggl[x-\frac12+\biggl(s+\frac14\biggr)^{1/2}\biggr] \biggl[\frac{\Gamma(x+1)}{\Gamma(x+s)}\biggr]^{1/(s-1)} \quad\text{and}\quad \biggl(x+\frac{s}2\biggr)\biggl[\frac{\Gamma(x+1)}{\Gamma(x+s)}\biggr]^{1/(s-1)}.
\end{equation}
\end{rem}

\section{Some completely monotonic functions involving ratios of two gamma or ${q}$-gamma functions}
\label{sec-cmf-ismail}

In this section, we look back and analyse several complete monotonicity of functions involving ratios of two gamma or $q$-gamma functions.

\subsection{Ismail-Lorch-Muldoon's monotonicity results}
Motivated by work on inequalities for the ratio of two gamma functions in~\cite{kershaw, laforgia-mc-1984, lorch-ultra} and~\cite[p.~155]{J.Wimp}, M.~E.~H. Ismail, L. Lorch and M.~E. Muldoon pointed out at the beginning of~\cite{Ismail-Lorch-Muldoon} that simple monotonicity of the ratio of two gamma functions are useful.
\par
In~\cite[pp.~118\nobreakdash--119]{Oliver}, the asymptotic formula
\begin{equation}\label{Oliver-asymp-formula-gamma-ratio}
z^{b-a}\frac{\Gamma(z+a)}{\Gamma(z+b)}\sim1+\frac{(a-b)(a+b-1)}{2z} +\dotsm
\end{equation}
as $z\to\infty$ along any curve joining $z=0$ and $z=\infty$ is listed, where $z\ne-a,-a-1,\dotsc$ and $z\ne-b,-b-1,\dotsc$. Suggested by it, the following complete monotonicity was proved in~\cite[Theorem~2.4]{Ismail-Lorch-Muldoon}: Let $a>b\ge0$, $a+b\ge1$ and
\begin{equation}\label{i-l-m-f-2-2}
h(x)=\ln\biggl[x^{a-b}\frac{\Gamma(x+b)}{\Gamma(x+a)}\biggr].
\end{equation}
Then both $h'(x)$ and
\begin{equation}\label{i-l-m-f-2}
x^{b-a}\frac{\Gamma(x+a)}{\Gamma(x+b)}
\end{equation}
are completely monotonic on $(0,\infty)$; The results fail when $a+b<1$ replaces $a+b\ge1$ in the hypotheses.
\par
Meanwhile, the following $q$-analogue of~\cite[Theorem~2.4]{Ismail-Lorch-Muldoon} was also provided in~\cite[Theorem~2.5]{Ismail-Lorch-Muldoon}: Let $a>b\ge0$, $a+b\ge1$, $q>0$, $q\ne1$ and
\begin{equation}
h_q(x)=\ln\biggl[|1-q^x|^{a-b}\frac{\Gamma_q(x+b)}{\Gamma_q(x+a)}\biggr].
\end{equation}
Then $h_q'(x)$ is completely monotonic on $(0,\infty)$; So is the function
\begin{equation}\label{i-l-m-f-2-q}
|1-q^x|^{b-a}\frac{\Gamma_q(x+a)}{\Gamma_q(x+b)};
\end{equation}
The result fails if $a+b<1$.

\begin{rem}
The proof of~\cite[Theorem~2.4]{Ismail-Lorch-Muldoon} can be outlined as follows: Using the integral representation
\begin{equation}\label{gauss-formula-psi}
\psi(z)=-\gamma+\int_0^\infty\frac{e^{-t}-e^{-tz}}{1-e^{-t}}\td t
\end{equation}
for $\re z>0$ yields
\begin{equation}\label{h'(x)-int}
h'(x)=\int_0^\infty\biggl[\frac{e^{-at}-e^{-bt}}{1-e^{-t}}+a-b\biggr]e^{-xt}\td t.
\end{equation}
It was established in~\cite[Lemma~4.1]{Ismail-Lorch-Muldoon} that if $0\le b<a$, $a+b\ge1$ and $b^2+(a-1)^2\ne0$, then
\begin{equation}\label{ILM-ineq-exp}
\frac{w^b-w^a}{1-w}<a-b,\quad 0<w<1;
\end{equation}
The result fails if the condition $a+b\ge1$ is replaced by $a+b<1$. Combining~\eqref{h'(x)-int} and~\eqref{ILM-ineq-exp} with Theorem~\ref{p.83-bochner} results in~\cite[Theorem~2.4]{Ismail-Lorch-Muldoon}.
\par
The proof of~\cite[Theorem~2.5]{Ismail-Lorch-Muldoon} was finished by using the formula~\eqref{q-gamma-1.4}, the inequality~\eqref{ILM-ineq-exp}, Theorem~\ref{p.161-widder} and Theorem~\ref{p.83-bochner}.
\end{rem}

\begin{rem}
It is noted that~\cite[Theorem~2.4 and Theorem~2.5]{Ismail-Lorch-Muldoon} mentioned above can be restated using the terminology ``logarithmically completely monotonic function'' as follows: The functions defined by~\eqref{i-l-m-f-2} and~\eqref{i-l-m-f-2-q} are logarithmically completely monotonic on $(0,\infty)$ if and only if $a+b\ge1$ for $a>b\ge0$, $q>0$ and $q\ne1$.
\end{rem}

\subsection{Bustoz-Ismail's monotonicity results}\label{Bustoz-Ismail-sec}
In~\cite{Bustoz-and-Ismail}, it was noticed that inequalities like~\eqref{theta-l-u-b} are ``immediate consequences of the complete monotonicity of certain functions. Indeed, one should investigate monotonicity properties of functions involving quotients of gamma functions and as a by-product derive inequalities of the aforementioned type. This approach is simpler and yields more general results.''
\par
In~\cite{Bustoz-and-Ismail}, it was revealed that
\begin{enumerate}
\item
\cite[Theorem~1]{Bustoz-and-Ismail}: the function
\begin{equation}\label{ismail-f(x)-thm1}
\frac1{(x+c)^{1/2}}\cdot\frac{\Gamma(x+1)}{\Gamma(x+1/2)},\quad x>\max\biggl\{-\frac12,-c\biggr\}
\end{equation}
is completely monotonic on $(-c,\infty)$ if $c\le\frac14$, so is the reciprocal of~\eqref{ismail-f(x)-thm1} on $\bigl[-\frac12,\infty\bigr)$ if $c\ge\frac12$;
\item
\cite[Theorem~3]{Bustoz-and-Ismail}: the function
\begin{equation}\label{ismail-g(x)-thm1}
(x+c)^{a-b}\frac{\Gamma(x+b)}{\Gamma(x+a)}
\end{equation}
for $1\ge b-a>0$ is completely monotonic on the interval $(\max\{-a,-c\},\infty)$ if $c\le\frac{a+b-1}2$, so is the reciprocal of~\eqref{ismail-g(x)-thm1} on $(\max\{-a,-c\},\infty)$ if $c\ge a$;
\item
\cite[Theorem~7]{Bustoz-and-Ismail}: the function
\begin{equation}\label{ismail-h(x)-thm2}
\frac{\Gamma(x+1)}{\Gamma(x+s)}\biggl(x+\frac{s}2\biggr)^{s-1}
\end{equation}
for $0\le s\le1$ is completely monotonic on $(0,\infty)$; when $0<s<1$, it satisfies $(-1)^nf^{(n)}(x)>0$ for $x>0$;
\item
\cite[Theorem~8]{Bustoz-and-Ismail}: the function
\begin{equation}\label{thm8-ismail-one}
\Biggl(x-\frac12+\sqrt{s+\frac14}\,\Biggr)^{1-s}\frac{\Gamma(x+s)}{\Gamma(x+1)}
\end{equation}
for $0<s<1$ is strictly decreasing on $(0,\infty)$.
\end{enumerate}

\begin{rem}\label{lemma2-bus-ism}
A special case of Theorem~\ref{p.83-bochner} says that the function $\exp(-h(x))$ is completely monotonic on an interval $I$ if $h'(x)$ is completely monotonic on $I$. This was iterated as \cite[Lemma~2.1]{Bustoz-and-Ismail}.
\par
In \cite[p.~15 and p.~20]{er}, the following integral representation was listed: For $\re z>0$,
\begin{equation}\label{psi-frac12}
\psi\biggl(\frac12+\frac{z}2\biggr)-\psi\biggl(\frac{z}2\biggr)
=2\int_0^\infty\frac{e^{-zt}}{1+e^{-t}}\td t.
\end{equation}
The formula~\eqref{psi-frac12} and \cite[Lemma~2.1]{Bustoz-and-Ismail} are basic tools of the proof of \cite[Theorem~1]{Bustoz-and-Ismail}.
\end{rem}

\begin{rem}
The basic tools of the proof of~\cite[Theorem~3]{Bustoz-and-Ismail} include~\cite[Lemma~2.1]{Bustoz-and-Ismail} mentioned in Remark~\ref{lemma2-bus-ism}, the formula~\eqref{gauss-formula-psi}, and the non-negativeness of the function
\begin{equation}\label{omega-ismail-nonneg-1}
\omega(t)=2(b-a)\sinh\frac{t}2-2\sinh\frac{(b-a)t}2
\end{equation}
for $b>a$ and $t\ge0$ and the function
\begin{equation}\label{omega-ismail-nonneg-2}
(a-b)(1-e^{-t})+e^{(c-a)t}-e^{(c-b)t}
\end{equation}
for $b>a$, $c\ge a$ and $t\ge0$.
\end{rem}

\begin{rem}
The proof of the complete monotonicity of the function~\eqref{ismail-h(x)-thm2} in~\cite[Theorem~7]{Bustoz-and-Ismail} relies on the series representation
\begin{equation}\label{series-repr}
\psi(x)=-\gamma-\frac1x+\sum_{n=1}^\infty\biggl(\frac1n-\frac1{x+n}\biggr)
\end{equation}
in \cite[p.~15]{er}, the positivity of the function
\begin{equation}\label{sinh-sinh-ismail}
(1-s)\sinh t-\sinh[(1-s)t]
\end{equation}
on $(0,\infty)$ for $0<s<1$, and the above Theorem~\ref{p.83-bochner} applied to $f(x)=e^{-x}$, as mentioned in Remark~\ref{lemma2-bus-ism}.
\par
The proof of the decreasing monotonicity of the function~\eqref{thm8-ismail-one} just used the formula~\eqref{series-repr} and the conclusion stated in Remark~\ref{lemma2-bus-ism}.
\end{rem}

\begin{rem}
In fact, under corresponding assumptions, the functions~\eqref{ismail-f(x)-thm1}, \eqref{ismail-g(x)-thm1}, \eqref{ismail-h(x)-thm2} and their reciprocals had been proved in~\cite{Bustoz-and-Ismail} to be logarithmically completely monotonic.
\end{rem}

\subsection{Ismail-Muldoon's monotonicity results}
It was claimed in~\cite[p.~310]{Ismail-Muldoon-119} that ``Many inequalities for special functions follow from monotonicity properties. Often such inequalities are special cases of the complete monotonicity of related special functions. For example, an inequality of the form $f(x)\ge g(x)$ for $x\in[a,\infty)$ with equality if and only if $x=a$ may be a disguised form of the complete monotonicity of $\frac{g(\varphi(x))}{f(\varphi(x))}$ where $\phi(x)$ is a nondecreasing function on $(a,\infty)$ and $\frac{g(\varphi(a))}{f(\varphi(a))}=1$''.
\par
Among other things, suggested by \cite[Theorem~3]{Bustoz-and-Ismail} mentioned in the above section, the following complete monotonicity was presented in \cite[Theorem~2.5]{Ismail-Muldoon-119}: Let $a<b\le a+1$ and
\begin{equation}
g(x)=\biggl(\frac{1-q^{x+c}}{1-q}\biggr)^{a-b}\frac{\Gamma_q(x+b)}{\Gamma_q(x+a)}.
\end{equation}
Then $-[\ln g(x)]'$ is completely monotonic on $(-c,\infty)$ if $0\le c\le\frac{a+b-1}2$ and $[\ln g(x)]'$ is completely monotonic on $(-a,\infty)$ if $c\ge a\ge0$; Neither is completely monotonic for $\frac{a+b-1}2<c<a$.
\par
As a supplement of \cite[Theorem~2.5]{Ismail-Muldoon-119}, it was proved separately in \cite[Theorem~2.6]{Ismail-Muldoon-119} that the first derivative of the function
\begin{equation}
\ln\biggl[\biggl(\frac{1-q^x}{1-q}\biggr)^a\frac{\Gamma_q(x)}{\Gamma_q(x+a)}\biggr],\quad 0<q<1
\end{equation}
is completely monotonic on $(0,\infty)$ for $a\ge1$.

\begin{rem}
The proof of \cite[Theorem~2.5]{Ismail-Muldoon-119} depends on deriving
\begin{equation}
\frac{\td{}}{\td x}\ln g(x)=-\int_0^\infty e^{-xt} \biggl[\frac{e^{-bt}-e^{-at}}{1-e^{-t}}+(b-a)e^{-ct}\biggr]\td\gamma_q(t)
\end{equation}
and \cite[Lemma~1.2]{Ismail-Muldoon-119}: Let $0<\alpha<1$. Then
\begin{equation}\label{Lemma1.2-Ismail-Muldoon-119}
\alpha e^{(\alpha-1)t}<\frac{\sinh(\alpha t)}{\sinh t}<\alpha,\quad t>0.
\end{equation}
The inequalities become equalities when $\alpha=1$ and they are reversed when $\alpha>1$.
\par
The proof of \cite[Theorem~2.6]{Ismail-Muldoon-119} is similar to that of \cite[Theorem~2.5]{Ismail-Muldoon-119}.
\end{rem}

\begin{rem}
It is clear that that Theorem~2.5 and Theorem~2.6 in~\cite{Ismail-Muldoon-119} mentioned above can be rewritten using the phrase ``logarithmically completely monotonic function''.
\end{rem}

\begin{rem}
From~\cite[Theorem~2.5]{Ismail-Muldoon-119}, the following inequality was derived in~\cite[Theorem~3.3]{Ismail-Muldoon-119}: For $0<q\le1$, the inequality
\begin{equation}\label{i-m-q-gamma-ineq}
\frac{\Gamma_q(x+1)}{\Gamma_q(x+s)}>\biggl(\frac{1-q^{x+s/2}}{1-q}\biggr)^{1-s}
\end{equation}
holds for $0<s<1$ and $x>-\frac{s}2$. In \cite{Alzer-Math-Nachr-2001}, it was pointed out that the inequality
\begin{equation}\label{alzer-q-gamma-ineq}
\frac{\Gamma_q(x+1)}{\Gamma_q(x+s)}<\biggl(\frac{1-q^{x+s}}{1-q}\biggr)^{1-s},\quad s\in(0,1)
\end{equation}
is also valid for $x>-s$. As refinements of~\eqref{i-m-q-gamma-ineq} and~\eqref{alzer-q-gamma-ineq}, the following double inequality was presented in~\cite[Theorem~3.1]{Alzer-Math-Nachr-2001}: For real numbers $0<q\ne1$ and $s\in(0,1)$, the double inequality
\begin{equation}
\biggl[\frac{1-q^{x+\alpha(q,s)}}{1-q}\biggr]^{1-s} <\frac{\Gamma_q(x+1)}{\Gamma_q(x+s)} <\biggl[\frac{1-q^{x+\beta(q,s)}}{1-q}\biggr]^{1-s},\quad x>0
\end{equation}
holds with the best possible values
\begin{equation}
\alpha(q,s)=\begin{cases}
\dfrac{\ln[(q^s-q)/(1-s)(1-q)]}{\ln q},&0<q<1\\
\dfrac{s}2,&q>1
\end{cases}
\end{equation}
and
\begin{equation}
\beta(q,s)=\frac{\ln\bigl\{1-(1-q)[\Gamma_q(s)]^{1/(s-1)}\bigr\}}{\ln q}.
\end{equation}
As a direct consequence, it was derived in~\cite[Corollary~3.2]{Alzer-Math-Nachr-2001} that the inequality
\begin{equation}\label{lug-egp-alzer-ineq}
[x+a(s)]^{1-s}\le\frac{\Gamma(x+1)}{\Gamma(x+s)}\le [x+b(s)]^{1-s}
\end{equation}
holds for $s\in(0,1)$ and $x\ge0$ with the best possible values $a(s)=\frac{s}2$ and $b(s)=[\Gamma(s)]^{1/(s-1)}$.
\par
The inequality~\eqref{lug-egp-alzer-ineq} was ever claimed in~\cite[p.~248]{Lazarevic} with a wrong proof. It was also generalized and extended in \cite[Theorem~3]{egp}.
\end{rem}

\section{Some logarithmically completely monotonic functions involving ratios of two gamma or ${q}$-gamma functions}

In this section, we look back and analyse necessary and sufficient conditions for functions involving ratios of two gamma or $q$-gamma functions to be logarithmically completely monotonic.

\subsection{Some properties of a function involving exponential functions}
For real numbers $\alpha$ and $\beta$ with $\alpha\ne\beta$ and $(\alpha,\beta)\not\in\{(0,1),(1,0)\}$, let
\begin{equation}\label{q-dfn}
q_{\alpha,\beta}(t)=
\begin{cases}
\dfrac{e^{-\alpha t}-e^{-\beta t}}{1-e^{-t}},&t\ne0;\\
\beta-\alpha,&t=0.
\end{cases}
\end{equation}
As seen in Section~\ref{sec-cmf-ismail}, it is easy to have an idea that the function $q_{\alpha,\beta}(t)$ or its variations play indispensable roles in the proofs of \cite[Theorem~3]{Bustoz-and-Ismail}, \cite[Theorem~7]{Bustoz-and-Ismail}, \cite[Theorem~2.4]{Ismail-Lorch-Muldoon}, \cite[Theorem~2.5]{Ismail-Lorch-Muldoon}, \cite[Theorem~2.5]{Ismail-Muldoon-119} and \cite[Theorem~2.6]{Ismail-Muldoon-119}.
\par
In order to bound ratios of two gamma or $q$-gamma functions, necessary and sufficient conditions for $q_{\alpha,\beta}(t)$ to be either monotonic or logarithmically convex have been investigated in~\cite{mon-element-exp-final.tex, mon-element-exp.tex-rgmia, comp-mon-element-exp.tex, notes-best-new-proof.tex, notes-best.tex-mia, notes-best.tex-rgmia}.

\begin{prop}[\cite{mon-element-exp-final.tex}]\label{q-mon-lem-2}
Let $t$, $\alpha$ and $\beta$ with $\alpha\ne\beta$ and $(\alpha,\beta)\not\in\{(0,1),(1,0)\}$ be real numbers. Then
\begin{enumerate}
\item
the function $q_{\alpha,\beta}(t)$ increases on $(0,\infty)$ if and only
if $(\beta-\alpha)(1-\alpha-\beta)\ge0$ and $(\beta-\alpha) (|\alpha-\beta|
-\alpha-\beta)\ge0$;
\item
the function $q_{\alpha,\beta}(t)$ decreases on $(0,\infty)$ if and only
if $(\beta-\alpha)(1-\alpha-\beta)\le0$ and $(\beta-\alpha) (|\alpha-\beta|
-\alpha-\beta)\le0$;
\item
the function $q_{\alpha,\beta}(t)$ increases on $(-\infty,0)$ if and only
if $(\beta-\alpha)(1-\alpha-\beta)\ge0$ and $(\beta-\alpha) (2-|\alpha-\beta|
-\alpha-\beta)\ge0$;
\item
the function $q_{\alpha,\beta}(t)$ decreases on $(-\infty,0)$ if and only
if $(\beta-\alpha)(1-\alpha-\beta)\le0$ and $(\beta-\alpha) (2-|\alpha-\beta|
-\alpha-\beta)\le0$;
\item
the function $q_{\alpha,\beta}(t)$ increases on $(-\infty,\infty)$ if and
only if $(\beta-\alpha) (|\alpha-\beta| -\alpha-\beta)\ge0$ and
$(\beta-\alpha) (2-|\alpha-\beta| -\alpha-\beta)\ge0$;
\item
the function $q_{\alpha,\beta}(t)$ decreases on $(-\infty,\infty)$ if and
only if $(\beta-\alpha) (|\alpha-\beta| -\alpha-\beta)\le0$ and
$(\beta-\alpha) (2-|\alpha-\beta| -\alpha-\beta)\le0$.
\end{enumerate}
\end{prop}

\begin{prop}[{\cite{mon-element-exp-final.tex}, \cite[Lemma~1]{notes-best.tex-mia} and~\cite[Lemma~1]{notes-best.tex-rgmia}}]\label{q-log-conv-thm}
The function $q_{\alpha,\beta}(t)$ on $(-\infty,\infty)$ is logarithmically
convex if $\beta-\alpha>1$ and logarithmically concave if $0<\beta-\alpha<1$.
\end{prop}

\begin{prop}[{\cite[Theorem~1.1]{comp-mon-element-exp.tex}}]
If $1>\beta-\alpha>0$, then $q_{\alpha,\beta}(u)$ is $3$-log-convex on $(0,\infty)$ and $3$-log-concave on $(-\infty,0)$; If $\beta-\alpha>1$, then $q_{\alpha,\beta}(u)$ is $3$-log-concave on $(0,\infty)$ and $3$-log-convex on $(-\infty,0)$.
\end{prop}

\begin{prop}[{\cite[Lemma~3]{notes-best-new-proof.tex}}]
Let $\lambda\in\mathbb{R}$. If $\beta-\alpha>1$, then the function $q_{\alpha,\beta}(t) q_{\alpha,\beta}(\lambda-t)$ is increasing on $\bigl(\frac\lambda2,\infty\bigr)$ and decreasing on $\bigl(-\infty, \frac\lambda2\bigr)$; if $0<\beta-\alpha<1$, it is decreasing on $\bigl(\frac\lambda2, \infty\bigr)$ and increasing on $\bigl(-\infty, \frac\lambda2\bigr)$.
\end{prop}

\begin{rem}
By noticing that the function $q_{\alpha,\beta}(t)$ can be rewritten as
\begin{equation}\label{rewr-f}
q_{\alpha,\beta}(t)=\frac{\sinh[(\beta-\alpha)t/2]}{\sinh(t/2)}\exp\frac{(1-\alpha-\beta)t}2,
\end{equation}
it is easy to see that the inequality~\eqref{ILM-ineq-exp}, the non-negativeness of the functions~\eqref{omega-ismail-nonneg-1} and~\eqref{omega-ismail-nonneg-2}, the positivity of the function~\eqref{sinh-sinh-ismail} and the inequality~\eqref{Lemma1.2-Ismail-Muldoon-119} are at all special cases of the monotonicity of the function $q_{\alpha,\beta}(t)$ on $(0,\infty)$ stated in Proposition~\ref{q-mon-lem-2}.
\end{rem}

\subsection{Necessary and sufficient conditions related to the ratio of two gamma functions}

In this section, we survey necessary and sufficient conditions for some functions involving the ratio of two gamma functions to be logarithmically completely monotonic.

\subsubsection{}
The logarithmically complete monotonicity of the function
\begin{equation}\label{h-def-sandor}
h_a(x)=\frac{(x+a)^{1-a}\Gamma(x+a)}{x\Gamma(x)}=\frac{(x+a)^{1-a}\Gamma(x+a)}{\Gamma(x+1)}
\end{equation}
for $x>0$ and $a>0$, the reciprocal of the first function in~\eqref{2.5-function} discussed in Remark~\ref{rem-2.1.1}, were considered in~\cite{sandor-gamma.tex-rgmia, sandor-gamma-JKMS.tex}.

\begin{thm}[\cite{sandor-gamma.tex-rgmia} and {\cite[Theorem~1.2]{sandor-gamma-JKMS.tex}}]\label{thm-sandor-qi-2}
The function $h_a(x)$ has the following properties:
\begin{enumerate}
\item
The function $h_a(x)$ is logarithmically completely monotonic on $(0,\infty)$ if $0<a<1$;
\item
The function $[h_a(x)]^{-1}$ is logarithmically completely monotonic on $(0,\infty)$ if $a>1$;
\item
For any $a>0$,
\begin{equation}
\lim_{x\to0^+}h_a(x)=\frac{\Gamma(a+1)}{a^a}
\quad \text{and}\quad \lim_{x\to\infty}h_a(x)=1.
\end{equation}
\end{enumerate}
\end{thm}

In order to obtain a refined upper bound in \eqref{wendel-inequal}, the logarithmically complete monotonicity of the function
\begin{equation}\label{f-def-sandor}
f_a(x)=\frac{\Gamma(x+a)}{x^a\Gamma(x)}
\end{equation}
for $x\in(0,\infty)$ and $a\in(0,\infty)$, the middle term in \eqref{wendel-inequal} or the reciprocal of the second function in~\eqref{2.5-function}, were considered in~\cite{sandor-gamma.tex-rgmia} and \cite[Theorem~1.3]{sandor-gamma-JKMS.tex}.

\begin{thm}[\cite{sandor-gamma.tex-rgmia} and {\cite[Theorem~1.3]{sandor-gamma-JKMS.tex}}] \label{thm-sandor-qi-final}
The function $f_a(x)$ has the following properties:
\begin{enumerate}
\item
The function $f_a(x)$ is logarithmically completely monotonic on $(0,\infty)$ and $\lim_{x\to0+}f_a(x)=\infty$ if $a>1$;
\item
The function $[f_a(x)]^{-1}$ is logarithmically completely monotonic on $(0,\infty)$ and $\lim_{x\to0+}f_a(x)=0$ if $0<a<1$;
\item
$\lim_{x\to\infty}f_a(x)=1$ for any $a\in(0,\infty)$.
\end{enumerate}
\end{thm}

As a straightforward consequence of combining Theorem~\ref{thm-sandor-qi-2} and Theorem~\ref{thm-sandor-qi-final}, the following refinement of the upper bound in the inequality~\eqref{wendel-inequal} is established.

\begin{thm}[\cite{sandor-gamma.tex-rgmia} and {\cite[Theorem~1.4]{sandor-gamma-JKMS.tex}}] \label{sandor-qi-inequal-ref}
Let $x\in(0,\infty)$. If $0<a<1$, then
\begin{multline}\label{combined-inequal}
\biggl(\frac{x}{x+a}\biggr)^{1-a} <\frac{\Gamma(x+a)}{x^a\Gamma(x)}\\*
<\begin{cases} \dfrac{\Gamma(a+1)}{a^a}\biggl(\dfrac{x}{x+a}\biggr)^{1-a}\le1,
&0<x\le\dfrac{ap(a)}{1-p(a)},
\\1,&\dfrac{ap(a)}{1-p(a)}<x<\infty,
\end{cases}
\end{multline}
where
\begin{equation}\label{p-def-sandor}
p(x)=\begin{cases}
\biggl[\dfrac{x^x}{\Gamma(x+1)}\biggr]^{1/(1-x)},&x\ne1,\\
e^{-\gamma},&x=1.
\end{cases}
\end{equation}
If $a>1$, the reversed inequality of~\eqref{combined-inequal} holds.
\end{thm}

\begin{rem}
The logarithmically complete monotonicity of the function~\eqref{p-def-sandor} and its generalized form were researched in~\cite{sandor-gamma.tex-rgmia}, \cite[Theorem~1.5]{sandor-gamma-JKMS.tex}, \cite[Theorem~1.4]{sandor-gamma-2-ITSF.tex} and~\cite[Theorem~1.4]{sandor-gamma-2-ITSF.tex-rgmia} respectively.
\end{rem}

\subsubsection{}
In~\cite[Theorem~1]{laj-7.pdf}, the following logarithmically complete monotonicity were established: The functions
\begin{equation}\label{laj-7-2-funct}
\frac{\Gamma(x+t)}{\Gamma(x+s)}\biggl(x+\frac{s+t-1}2\biggr)^{s-t}\quad \text{and}\quad \frac{\Gamma(x+s)}{\Gamma(x+t)}(x+s)^{t-s}
\end{equation}
for $0<s<t<s+1$ are logarithmically completely monotonic with respect to $x$ on $(-s,\infty)$.

\begin{rem}
We can not understand why the authors of the paper~\cite{laj-7.pdf} chose so special functions in~\eqref{laj-7-2-funct}. More accurately, we have no idea why the constants $\frac{s+t-1}2$ and $s$ were chosen in the polynomial factors of the functions listed in~\eqref{laj-7-2-funct}. Perhaps this can be interpreted by Theorem~\ref{unify-log-comp-thm} and Theorem~\ref{polygamma-divided} below.
\end{rem}

\subsubsection{}
For real numbers $a$, $b$ and $c$, denote $\rho=\min\{a,b,c\}$ and let
\begin{equation}\label{h-def-sandor-new}
H_{a,b;c}(x)=(x+c)^{b-a}\frac{\Gamma(x+a)}{\Gamma(x+b)}
\end{equation}
for $x\in(-\rho,\infty)$.
\par
By a recourse to the incomplete monotonicity of the function $q_{\alpha,\beta}(t)$ obtained in~\cite{mon-element-exp.tex-rgmia}, the following incomplete but correct conclusions about the logarithmically complete monotonicity of the function $H_{a,b;c}(x)$ were procured in~\cite{sandor-gamma-3.tex-jcam, sandor-gamma-3.tex-rgmia}.

\begin{thm}[{\cite[Theorem~1]{sandor-gamma-3.tex-jcam} and~\cite[Theorem~1]{sandor-gamma-3.tex-rgmia}}]\label{unify-log-comp-thm-orig}
Let $a$, $b$ and $c$ be real numbers and $\rho=\min\{a,b,c\}$. Then
\begin{enumerate}
\item
the function $H_{a,b;c}(x)$ is logarithmically completely monotonic on $(-\rho,\infty)$ if
\begin{equation*}
\begin{aligned}
(a,b;c)&\in\biggl\{(a,b;c):a+b\ge1,c\le b<c+\frac12\biggr\}
\cup\biggl\{(a,b;c):a>b\ge c+\frac12\biggr\}\\
&\quad\cup\{(a,b;c):2a+1\le a+b\le1,a<c\}
\cup\{(a,b;c):b-1\le a<b\le c\} \\
&\quad\setminus\{(a,b;c):a=c+1,b=c\},
\end{aligned}
\end{equation*}
\item
so is the function $[H_{a,b;c}(x)]^{-1}$ if
\begin{equation*}
\begin{aligned}
(a,b;c)&\in\biggl\{(a,b;c):a+b\ge1, c\le a<c+\frac12\biggr\}
\cup\biggl\{(a,b;c):b>a\ge c+\frac12\biggr\}\\
&\quad\cup\{(a,b;c):b<a\le c\}
\cup\{(a,b;c):b+1\le a,c\le a\le c+1\}\\
&\quad\cup\{(a,b;c):b+c+1\le a+b\le1\}\\
&\quad\setminus\{(a,b;c):a=c+1,b=c\}
\setminus\{(a,b;c):b=c+1,a=c\}.
\end{aligned}
\end{equation*}
\end{enumerate}
\end{thm}

\subsubsection{}
In~\cite[Theorem~1]{notes-best-simple-equiv.tex}, \cite[Theorem~1]{notes-best-simple-equiv.tex-RGMIA} and~\cite[Theorem~2]{notes-best-simple.tex-rgmia}, the function
\begin{equation}\label{differen-ineq}
\delta_{s,t}(x)=
\begin{cases}
\dfrac{\psi(x+t)-\psi(x+s)}{t-s}-\dfrac{2x+s+t+1}{2(x+s)(x+t)},&s\ne t\\[1em]
\psi'(x+s)-\dfrac1{x+s}-\dfrac1{2(x+s)^2},&s=t
\end{cases}
\end{equation}
for $|t-s|<1$ and $-\delta_{s,t}(x)$ for $|t-s|>1$ were proved to be completely monotonic on the interval $(-\min\{s,t\},\infty)$. By employing the formula~\eqref{gauss-formula-psi}, the monotonicity of $q_{\alpha,\beta}(t)$ on $(0,\infty)$ and the complete monotonicity of $\delta_{s,t}(x)$, necessary and sufficient conditions are presented for the function $H_{a,b;c}(x)$ to be logarithmically completely monotonic on $(-\rho,\infty)$ as follows.

\begin{thm}[\cite{sandor-gamma-3-note.tex-final, sandor-gamma-3-note.tex}]\label{unify-log-comp-thm}
Let $a$, $b$ and $c$ be real numbers and $\rho=\min\{a,b,c\}$. Then
\begin{enumerate}
\item
the function $H_{a,b;c}(x)$ is logarithmically completely monotonic on $(-\rho,\infty)$ if
and only if
\begin{equation}\label{d1-dfn-new}
\begin{split}
(a,b;c)\in D_1(a,b;c)&\triangleq\{(a,b;c):(b-a)(1-a-b+2c)\ge0\}\\
&\quad\cap\{(a,b;c):(b-a) (\vert a-b\vert-a-b+2c)\ge0\}\\
&\quad\setminus\{(a,b;c):a=c+1=b+1\}\\
&\quad\setminus\{(a,b;c):b=c+1=a+1\};
\end{split}
\end{equation}
\item
so is the function $H_{b,a;c}(x)$ on $(-\rho,\infty)$ if and only if
\begin{equation}\label{d2-dfn-new}
\begin{split}
(a,b;c)\in D_2(a,b;c)&\triangleq\{(a,b;c):(b-a)(1-a-b+2c)\le0\}\\
&\quad\cap\{(a,b;c):(b-a) (\vert a-b\vert-a-b+2c)\le0\}\\
&\quad\setminus\{(a,b;c):b=c+1=a+1\}\\
&\quad\setminus\{(a,b;c):a=c+1=b+1\}.
\end{split}
\end{equation}
\end{enumerate}
\end{thm}

\begin{rem}
The limit~\eqref{wendel-approx} implies that $\lim_{x\to\infty}H_{a,b;c}(x)=1$ is valid for all defined numbers $a,b,c$. Combining this with the logarithmically complete monotonicity of $H_{a,b;c}(x)$ yields that the inequality
\begin{equation}
H_{a,b;c}(x)>1
\end{equation}
holds if $(a,b;c)\in D_1(a,b;c)$ and reverses if $(a,b;c)\in D_2(a,b;c)$, that is, the inequality
\begin{equation}\label{wendel-gamma-ineq-orig}
x+\lambda<\biggl[\frac{\Gamma(x+a)}{\Gamma(x+b)}\biggr]^{1/(a-b)}<x+\mu,\quad b>a
\end{equation}
holds for $x\in(-a,\infty)$ if $\lambda\le\min\bigl\{a,\frac{a+b-1}2\bigr\}$ and $\mu\ge\max\bigl\{a,\frac{a+b-1}2\bigr\}$, which is equivalent to
\begin{equation}\label{wendel-gamma-ineq}
\min\biggl\{a,\frac{a+b-1}2\biggr\}<\biggl[\frac{\Gamma(a)}{\Gamma(b)}\biggr]^{1/(a-b)} <\max\biggl\{a,\frac{a+b-1}2\biggr\},\quad b>a>0.
\end{equation}
\par
It is noted that a special case $0<a<b<1$ of the inequality~\eqref{wendel-gamma-ineq-orig} was derived in~\cite{Chen-Oct-04-1051} from Elezovi\'c-Giordano-Pe\v{c}ari\'c's theorem (see \cite{egp, notes-best-new-proof.tex, notes-best.tex-mia, notes-best.tex-rgmia}). Moreover, by available of the inequality~\eqref{wendel-inequal} and others, the double inequalities
\begin{equation}
\frac{x+a}{x+b}(x+b)^{b-a}\le\frac{\Gamma(x+b)}{\Gamma(x+a)}\le(x+a)^{b-a},\quad x>0
\end{equation}
and
\begin{equation}
(x+a)e^{-\gamma/(x+a)}<\biggl[\frac{\Gamma(x+b)}{\Gamma(x+a)}\biggr]^{1/(b-a)} <(x+b)e^{-1/2(x+b)}, \quad x\ge1
\end{equation}
were proved in~\cite{Sandor-Oct-04-1052} to be valid for $0<a<b<1$.
\par
Maybe two references~\cite{Bencze-OQ1352, Modan-Oct-04-1055} are also useful and worth being mentioned.
\end{rem}

\begin{rem}
Since the complete monotonicity of the function~\eqref{differen-ineq} was not established and the main result in~\cite{mon-element-exp.tex-rgmia} about the monotonicity of the function $q_{\alpha,\beta}(t)$ is incomplete at that time, necessary conditions for the function~\eqref{h-def-sandor-new} to be logarithmically completely monotonic was not discovered in~\cite[Theorem~1]{sandor-gamma-3.tex-jcam} and~\cite[Theorem~1]{sandor-gamma-3.tex-rgmia} and the sufficient conditions in~\cite[Theorem~1]{sandor-gamma-3.tex-jcam} and~\cite[Theorem~1]{sandor-gamma-3.tex-rgmia} are imperfect.
\end{rem}

\begin{rem}
It is not difficult to see that all (complete) monotonicity on functions involving the ratio of two gamma functions, showed by Bustoz-Ismail in~\cite{Bustoz-and-Ismail} and Ismail-Lorch-Muldoon in~\cite{Ismail-Lorch-Muldoon} and related results in \cite{sandor-gamma-3.tex-jcam, sandor-gamma-3.tex-rgmia, sandor-gamma.tex-rgmia, sandor-gamma-JKMS.tex}, are special cases of the above Theorem~\ref{unify-log-comp-thm}.
\end{rem}

\subsubsection{}
From the above Theorem~\ref{unify-log-comp-thm}, the following double inequalities for divided differences of the psi and polygamma functions may be deduced immediately.

\begin{thm}[\cite{sandor-gamma-3-note.tex-final}]\label{polygamma-divided}
Let $b>a\ge0$ and $k\in\mathbb{N}$. Then the double inequality
\begin{equation}\label{n-s-ineq}
\frac{(k-1)!}{(x+\alpha)^k}<\frac{(-1)^{k-1} \bigl[\psi^{(k-1)}(x+b)-\psi^{(k-1)}(x+a)\bigr]}{b-a} <\frac{(k-1)!}{(x+\beta)^k}
\end{equation}
for $x\in(-\rho,\infty)$ holds if $\alpha\ge\max\bigl\{a,\frac{a+b-1}2\bigr\}$ and $0\le\beta\le\min\bigl\{a,\frac{a+b-1}2\bigr\}$.
\end{thm}

\begin{rem}
It is amazing that taking $b-a=1$ in~\eqref{n-s-ineq} leads to
\begin{equation}
\psi^{(k-1)}(x+a+1)-\psi^{(k-1)}(x+a)=(-1)^{k-1}\frac{(k-1)!}{(x+a)^k}
\end{equation}
for $a\ge0$, $x>0$ and $k\in\mathbb{N}$, which is equivalent to the recurrence formula
\begin{equation}\label{recurrence-formula}
\psi^{(n)}(z+1)-\psi^{(n)}(z)=(-1)^nn!z^{-n-1},\quad z>0,\quad n\ge0
\end{equation}
listed in~\cite[p.~260, 6.4.6]{abram}. For detailed information, see~\cite{roots-polygamma-eq.tex-ajmaa, roots-polygamma-eq.tex} and~\cite[Remark~8]{mon-element-exp-final.tex}.
\end{rem}

\begin{rem}
For more information on results of divided differences for the psi and polygamma functions, please refer to \cite{notes-best-simple-equiv.tex-RGMIA, notes-best-simple-equiv.tex, notes-best-simple-open.tex, simple-equiv.tex, notes-best-simple.tex-rgmia, simple-equiv-simple-rev.tex, Comp-Mon-Digamma-Trigamma-Divided.tex, AAM-Qi-09-PolyGamma.tex} and related references therein.
\end{rem}

\begin{rem}
It is worthwhile to note that some errors and defects appeared in~\cite{sandor-gamma-3.tex-jcam, sandor-gamma-3.tex-rgmia} have been corrected and consummated in~\cite{sandor-gamma-3-note.tex-final, sandor-gamma-3-note.tex}.
\end{rem}

\subsection{Necessary and sufficient conditions related to the ratio of two ${q}$-gamma functions}
The known results obtained by many mathematicians show that most of properties of the ratio of two gamma functions may be replanted to cases of the ratio of two $q$-gamma functions, as done in~\cite[Theorem~2.5]{Ismail-Lorch-Muldoon} and~\cite[Theorem~2.5 and Theorem~2.6]{Ismail-Muldoon-119} mentioned above.
\par
Let $a,b$ and $c$ be real numbers, $\rho=\min\{a,b,c\}$, and define
\begin{equation}\label{H{q;a,b;c}(x)}
H_{q;a,b;c}(x)=\biggl(\frac{1-q^{x+c}}{1-q}\biggr)^{a-b}\frac{\Gamma_q(x+b)}{\Gamma_q(x+a)}
\end{equation}
for $x\in(-\rho,\infty)$, where $\Gamma_q(x)$ is the $q$-gamma function defined by~\eqref{q-gamma-dfn} and~\eqref{q-gamma-dfn-q>1}.
\par
It is clear that the function~\eqref{H{q;a,b;c}(x)} is a $q$-analogue of the function~\eqref{h-def-sandor-new}.
\par
In virtue of the monotonicity of $q_{\alpha,\beta}(t)$ on $(0,\infty)$ and the formula~\eqref{q-gamma-1.5}, the following Theorem~\ref{q-gamma-ratio}, a $q$-analogue of Theorem~\ref{unify-log-comp-thm}, was procured.

\begin{thm}[\cite{mon-element-exp-final.tex}]\label{q-gamma-ratio}
Let $a$, $b$ and $c$ be real numbers and $\rho=\min\{a,b,c\}$. Then the function $H_{q;a,b;c}(x)$ is logarithmically completely monotonic on $(-\rho,\infty)$ if and only if $(a,b;c)\in D_2(a,b;c)$, so is the function $H_{q;b,a;c}(x)$ if and only if  $(a,b;c)\in D_1(a,b;c)$, where $D_1(a,b;c)$ and $D_2(a,b;c)$ are defined by~\eqref{d1-dfn-new} and~\eqref{d2-dfn-new} respectively.
\end{thm}

\begin{rem}
All complete monotonicity obtained in~\cite[Theorem~2.5]{Ismail-Lorch-Muldoon} and~\cite[Theorem~2.5 and Theorem~2.6]{Ismail-Muldoon-119} are special cases of Theorem~\ref{q-gamma-ratio}.
\end{rem}

Similar to Theorem~\ref{polygamma-divided}, the following double inequality of divided differences of the $q$-psi function $\psi_q(x)$ for $0<q<1$ may be derived from Theorem~\ref{q-gamma-ratio}.

\begin{thm}[\cite{mon-element-exp-final.tex}]\label{polygamma-divided-q-gamma}
Let $b>a\ge0$, $k\in\mathbb{N}$ and $0<q<1$. Then the inequality
\begin{equation}\label{n-s-ineq-g-gamma}
\frac{(-1)^{k-1}\bigl[\psi^{(k-1)}_q(x+b)-\psi^{(k-1)}_q(x+a)\bigr]}{b-a} <(-1)^{k-1}[\ln(1-q^{x+c})]^{(k)}
\end{equation}
for $x\in(-\rho,\infty)$ holds if $0\le c\le\min\bigl\{a,\frac{a+b-1}2\bigr\}$ and reverses if $c\ge\max\bigl\{a,\frac{a+b-1}2\bigr\}$. Consequently, the identity
\begin{equation}\label{n-s-ineq-g-gamma-equality}
\psi^{(k-1)}_q(x+1)-\psi^{(k-1)}_q(x) =[\ln(1-q^x)]^{(k)}
\end{equation}
holds for $x\in(0,\infty)$ and $k\in\mathbb{N}$.
\end{thm}

\begin{rem}
Since identities~\eqref{recurrence-formula} and~\eqref{n-s-ineq-g-gamma-equality} may be derived from inequalities~\eqref{n-s-ineq} and~\eqref{n-s-ineq-g-gamma}, we can regard inequalities~\eqref{n-s-ineq} and~\eqref{n-s-ineq-g-gamma} as generalizations of identities~\eqref{recurrence-formula} and~\eqref{n-s-ineq-g-gamma-equality}.
\end{rem}

\section{Logarithmically complete monotonicity for ratios of products of the gamma and ${q}$-gamma functions}

In this section, we would like to look back and analyse some (logarithmically) complete monotonicity of ratios of products of the gamma and $q$-gamma functions.
\par
Let $a_i$ and $b_i$ for $1\le i\le n$ be real numbers and $\rho_n=\min_{1\le i\le n}\{a_i,b_i\}$. For $x\in(-\rho_n,\infty)$, define
\begin{equation}
h_{\boldsymbol{a},\boldsymbol{b};n}(x)=\prod_{i=1}^n\frac{\Gamma(x+a_i)}{\Gamma(x+b_i)},
\end{equation}
where $\boldsymbol{a}$ and $\boldsymbol{b}$ denote $(a_1,a_2,\dotsc,a_n)$ and $(b_1,b_2,\dotsc,b_n)$ respectively.

\subsection{Complete monotonicity}
In \cite[Theorem~6]{Bustoz-and-Ismail}, by virtue of the formula~\eqref{gauss-formula-psi} and a special case of Theorem~\ref{p.83-bochner} mentioned in Remark~\ref{lemma2-bus-ism} above, the function
\begin{equation}\label{ismail-bustoz-ratio-gamma}
x\mapsto\frac{\Gamma(x)\Gamma(x+a+b)}{\Gamma(x+a)\Gamma(x+b)}
\end{equation}
for $a,b\ge0$, a special cases of $h_{\boldsymbol{a},\boldsymbol{b};n}(x)$ for $n=2$, was proved to be completely monotonic on $(0,\infty)$.
\par
In \cite[Theorem~10]{psi-alzer}, the function $h_{\boldsymbol{a},\boldsymbol{b};n}(x)$ was proved to be completely monotonic on $(0,\infty)$ provided that $0\le a_1\le\dotsm\le a_n$, $0\le b_1 \le \dotsm \le b_n$ and $\sum_{i=1}^ka_i\le\sum_{i=1}^kb_i$ for $1\le k\le n$. Its proof used the formula~\eqref{gauss-formula-psi}, a special case of Theorem~\ref{p.83-bochner} applied to $f(x)=e^{-x}$, and the following conclusion cited from \cite[p.~10]{marolk}: Let $a_i$ and $b_i$ for $i=1,\dotsc,n$ be real numbers such that $a_1\le\dotsm\le a_n$, $b_1\le\dotsm\le b_n$, and $\sum_{i=1}^ka_i\le\sum_{i=1}^kb_i$ for $k=1,\dotsc,n$. If the function $f$ is decreasing and convex on $\mathbb{R}$, then
\begin{equation}
\sum_{i=1}^nf(b_i)\le\sum_{i=1}^nf(a_i).
\end{equation}
\par
In~\cite[Theorem~4.1]{Ismail-Muldoon-119}, the functions
\begin{equation}
-\frac{\td}{\td x}\ln\frac{\Gamma_q(x+a_1) \Gamma_q(x+a_2)\dotsm\Gamma_q(x+a_n)}{[\Gamma(x+\bar{a})]^n}
\end{equation}
and
\begin{equation}\label{no-mean-gamma}
\frac{\td}{\td x}\ln\frac{\Gamma_q(x+a_1) \Gamma_q(x+a_2)\dotsm\Gamma_q(x+a_n)} {[\Gamma_q(x)]^{n-1}\Gamma_q(x+a_1+a_2+\dotsm+a_n)}
\end{equation}
were proved to be completely monotonic on $(0,\infty)$, where $a_1,\dotsc,a_n$ are positive numbers, $n\bar{a}=a_1+\dotsm+a_n$, and $0<q\le1$.
\par
In~\cite{Malig-Pecaric-Persson-95}, the function
\begin{equation}\label{mal-pec-pers-f}
x\mapsto\frac{[\Gamma(x)]^{n-1}\Gamma\bigl(x+\sum_{i=1}^na_i\bigr)}{\prod_{i=1}^n\Gamma(x+a_i)}
\end{equation}
for $a_i>0$ and $i=1,\dotsc,n$ was found to be decreasing on $(0,\infty)$.
\par
Motivated by the decreasing monotonic property of the function~\eqref{mal-pec-pers-f}, H.~Alzer proved in \cite[Theorem~11]{psi-alzer} that the function
\begin{equation}
x\mapsto\frac{[\Gamma(x)]^\alpha\Gamma\bigl(x+\sum_{i=1}^na_i\bigr)}{\prod_{i=1}^n\Gamma(x+a_i)}
\end{equation}
is completely monotonic on $(0,\infty)$ if and only if $\alpha=n-1$.

\begin{rem}
It is clear that the decreasingly monotonic property of the function~\eqref{mal-pec-pers-f} is just the special case $q\to1^-$ of the complete monotonicity of the function~\eqref{no-mean-gamma}. Therefore, it seems that the authors of the papers \cite{psi-alzer, Malig-Pecaric-Persson-95} were not aware of the results in~\cite[Theorem~4.1]{Ismail-Muldoon-119}.
\end{rem}

\begin{rem}
The complete monotonicity mentioned just now are indeed logarithmically completely monotonic ones.
\end{rem}

\subsection{Logarithmically complete monotonicity}

Let $S_n$ be the symmetric group over $n$ symbols, $a_1,a_2,\dotsc,a_n$. Let $O_n$ and $E_n$ be the sets of odd and even permutations over $n$ symbols, respectively. For $a_1>a_2>\dotsm>a_n>0$, define
\begin{equation}
F(x)=\frac{\prod_{\sigma\in E_n}\Gamma\bigl(x+a_{\sigma(2)} +2a_{\sigma(3)}+\dotsm+(n-1)a_{\sigma(n)}\bigr)} {\prod_{\sigma\in O_n}\Gamma\bigl(x+a_{\sigma(2)} +2a_{\sigma(3)}+\dotsm+(n-1)a_{\sigma(n)}\bigr)}.
\end{equation}
It was proved in \cite[Theorem~1.1]{grin-ismail} that the function $F(x-a_2-2a_3-\dotsm-(n-1)a_n)$ is logarithmically completely monotonic on $(0,\infty)$.
\par
In~\cite[Theorem~1.2]{grin-ismail}, it was presented that the functions
\begin{equation}\label{fn-dfn}
F_n(x)=\frac{\Gamma(x)\prod_{k=1}^{[n/2]}\Bigl[\prod_{m\in P_{n,2k}}
\Gamma\Bigl(x+\sum_{j=1}^{2k}a_{m_j}\Bigr)\Bigr]}
{\prod_{k=1}^{[(n+1)/2]}\Bigl[\prod_{m\in P_{n,2k-1}}
\Gamma\Bigl(x+\sum_{j=1}^{2k-1}a_{m_j}\Bigr)\Bigr]}
\end{equation}
for any $a_k>0$ and $k\in\mathbb{N}$ are logarithmically completely monotonic on $(0,\infty)$ and that any product of functions of the type~\eqref{fn-dfn} with different parameters $a_k$ is logarithmically completely monotonic on $(0,\infty)$ as well, where $P_{n,k}$ for $1\le k\le n$ is the set of all vectors $\boldsymbol{m}=(m_1,\dotsc,m_k)$ whose components are natural numbers such that $1\le m_\nu<m_\mu\le n$ for $1\le\nu<\mu\le k$ and $P_{n,0}$ is the empty set.

\begin{rem}
The above Theorem~1.2 is more general than Theorem~1.1. The case $n=2$ in Theorem~1.2 corresponds to the complete monotonicity of the function~\eqref{ismail-bustoz-ratio-gamma} obtained in~\cite[Theorem~6]{Bustoz-and-Ismail}.
\end{rem}

In~\cite[Theorem~3.2]{grin-ismail}, it was showed that if
\begin{equation}
F_q(x)=\frac{\prod_{\sigma\in E_n}\Gamma_q\bigl(x+a_{\sigma(2)} +2a_{\sigma(3)}+\dotsm+(n-1)a_{\sigma(n)}\bigr)} {\prod_{\sigma\in O_n}\Gamma_q\bigl(x+a_{\sigma(2)} +2a_{\sigma(3)}+\dotsm+(n-1)a_{\sigma(n)}\bigr)}
\end{equation}
for $a_1>a_2>\dotsm>a_n>0$, then $F_q(x-a_2-2a_3-\dotsm-(n-1)a_n)$ is a logarithmically completely monotonic function of $x$ on $(0,\infty)$.
\par
In~\cite[Theorem~3.3]{grin-ismail}, it was stated that the functions
\begin{equation}\label{fn-dfn-2}
F_{n,q}(x)=\frac{\Gamma_q(x)\prod_{k=1}^{[n/2]}\Bigl[\prod_{m\in P_{n,2k}} \Gamma_q\Bigl(x+\sum_{j=1}^{2k}a_{m_j}\Bigr)\Bigr]}
{\prod_{k=1}^{[(n+1)/2]}\Bigl[\prod_{m\in P_{n,2k-1}} \Gamma_q\Bigl(x+\sum_{j=1}^{2k-1}a_{m_j}\Bigr)\Bigr]}
\end{equation}
for any $a_k>0$ with $k=1,\dotsc,n$ are logarithmically completely monotonic on $(0,\infty)$, so is any product of functions~\eqref{fn-dfn-2} with different parameters $a_k$.

\begin{rem}
It is obvious that~\cite[Theorem~3.2 and Theorem~3.3]{grin-ismail} are $q$-analogues of~\cite[Theorem~1.1 and Theorem~1.2]{grin-ismail}.
\end{rem}

\subsection{Some recent conclusions}

By a recourse to the monotonicity of $q_{\alpha,\beta}(t)$ on $(0,\infty)$, the following sufficient conditions for the function $h_{\boldsymbol{a},\boldsymbol{b};n}(x)$ to be logarithmically completely monotonic on $(0,\infty)$ are devised.

\begin{thm}[\cite{mon-element-exp-final.tex}]\label{products-ratio-thm8}
If
\begin{equation}\label{c-1-cond}
(b_i-a_i)(1-a_i-b_i)\ge0\quad\text{and}\quad (b_i-a_i) (|a_i-b_i|-a_i-b_i)\ge0
\end{equation}
hold for $1\le i\le n$ and
\begin{equation}\label{c-3-cond}
\sum_{i=1}^nb_i\ge\sum_{i=1}^na_i,
\end{equation}
then the function $h_{\boldsymbol{a},\boldsymbol{b};n}(x)$ is logarithmically completely monotonic on $(-\rho_n,\infty)$. If inequalities in~\eqref{c-1-cond} and~\eqref{c-3-cond} are reversed, then the function $h_{\boldsymbol{b},\boldsymbol{a};n}(x)$ is logarithmically completely monotonic on $(-\rho_n,\infty)$.
\end{thm}

The $q$-analogue of Theorem~\ref{products-ratio-thm8} is as follows.

\begin{thm}[\cite{mon-element-exp-final.tex}]\label{products-ratio-thm9}
Let $a_i$ and $b_i$ for $1\le i\le n$ be real and $\rho_n=\min_{1\le i\le n}\{a_i,b_i\}$. For $x\in(-\rho_n,\infty)$, define
\begin{equation}
h_{q;\boldsymbol{a},\boldsymbol{b};n}(x)=\prod_{i=1}^n\frac{\Gamma_q(x+a_i)}{\Gamma_q(x+b_i)}
\end{equation}
for $0<q<1$, where $\boldsymbol{a}$ and $\boldsymbol{b}$ stand for $(a_1,a_2,\dotsc,a_n)$ and $(b_1,b_2,\dotsc,b_n)$ respectively. If inequalities in~\eqref{c-1-cond} and~\eqref{c-3-cond} hold, then the function $h_{q;\boldsymbol{a},\boldsymbol{b};n}(x)$ is logarithmically completely monotonic on $(-\rho_n,\infty)$. If inequalities in~\eqref{c-1-cond} and~\eqref{c-3-cond} are reversed, then the function $h_{q;\boldsymbol{b},\boldsymbol{a};n}(x)$ is logarithmically completely monotonic on $(-\rho_n,\infty)$.
\end{thm}

\subsection*{Acknowledgements}
This article was ever reported on 9 October 2008 as a talk in the seminar held at the RGMIA, School of Engineering and Science, Victoria University, Australia, while the author was visiting the RGMIA between March 2008 and February 2009 by the grant from the China Scholarship Council. The author expresses many thanks to Professors Pietro Cerone and Server S.~Dragomir and other local colleagues at Victoria University for their invitation and hospitality throughout this period.

\end{document}